\documentclass[12pt]{amsart}

\newcommand{\CP}{{\mathbb {CP}}}

\newcommand{\R}{{\mathbb R}}
\newcommand{\C}{{\mathbb C}}

\begin{document}

\author{V.A.~Vassiliev}
\title{A few problems on monodromy and discriminants}
\date{}
\address{Steklov Mathematical Institute of Russian Academy of Sciences; National
Research Institute --- Higher School of Economics} \email{vva@mi.ras.ru}
\thanks{
Supported by Program ``Leading scientific schools'', grant No. NSh-5138.2014.1
and RFBR grant 13-01-00383}

\maketitle

\section{Explicit obstructions to the Lyashko--Looijenga covering (and its real
analogs) for non-simple singularities} \label{llob}

The so-called Lyashko--Looijenga covering (see \cite{looj}, \cite{Applied}) is
a strong tool for constructing (or proving the existence of) the perturbations
of simple singularities with prescribed topological properties, such as
singularity types of different critical points, or intersection matrices of
vanishing cycles, see e.g. \cite{Lya}. The real version of this tool allows one
to construct and enumerate all topologically different Morsifications of real
simple singularities, see \cite{looj2}, \cite{Chi}, \cite{alg}.

A large amount of these options is preserved for non-simple singularities, see
\cite{alg}, \cite{Applied}. In particular, this method has predicted the
existence of many Morsifications with prescribed properties and indicated their
topological characteristics, so that it was easy to give a strict construction
of these Morsifications. However, in this case this method is rater
experimental or heuristic, without clear guaranties that all perturbations
found by it actually do exist. Therefore it is important to fix the
restrictions of this method. Here are several explicit problems.

\subsection{Complex version}

Let $f: (\C^n,0) \to (\C^1,0)$ be an isolated holomorphic function singularity,
$\mu$ its Milnor number, $F(x,\lambda): (\C^n \times \C^\mu,0) \to (\C^1,0)$
the miniversal deformation of $f$, and $\Sigma \subset \C^\mu$ the complete
bifurcation set of functions of this deformation, i.e. the set of values of the
parameter $\lambda \in \C^\mu$ such that the corresponding function $f_\lambda
\equiv F(\cdot, \lambda)$ has less than $\mu$ different critical values at
critical points close to the origin. The {\em Lyashko--Looijenga map} sends any
point $\lambda$ from a small neighborhood $B_\varepsilon$ of the origin in
$\C^\mu$ to the unordered collection of critical values of the function
$f_\lambda$ at points close to $0 \in \C^n$ (or, which is equivalent but
sometimes more convenient, to the set of values of basic symmetric polynomials
of these critical values). If the singularity of $f$ is simple then the
restriction of this map to $B_\varepsilon \setminus \Sigma$ defines a local
covering over the configuration space $B(D,\mu)$ of all subsets of cardinality
$\mu$ in a very small (even with respect to $\varepsilon$) neighborhood $D$ of
the origin in $\C$, see \cite{looj}. In particular, any element $ \alpha \in
\pi_1(B(D,\mu))$ can be realised by a loop which can be lifted to a path in
$B_\varepsilon \setminus \Sigma$ covering this loop.

For non-simple singularities this is not more the case. As previously, the
Lyashko--Looijenga map is submersive (and hence locally bijective) everywhere
in $B_\varepsilon \setminus \Sigma$ (this follows from the very notion of
versality). However, a sufficiently complicated path in $B(D, \mu)$, being
lifted into $\C^\mu \setminus \Sigma$ in accordance with this local
bijectivity, can run out from the neighborhood of the origin in $\C^\mu$. This
is related with the fact that for non-simple singularities the
Lyashko--Looijenga map is not proper: the preimage of the collection $(0,
\dots, 0)$ is the entire (positive-dimensional) $\mu=const$ stratum\footnote{A
weaker substitute for the Lyashko--Looijenga covering theorem holds in the case
of parabolic singularities, if one writes the versal deformation in the
canonical monomial form and allows large travellings in the space $\C^\mu$, see
\cite{Jaw}. For more complicated singularities the situation is even worse.}.
\medskip

\noindent {\bf Problem \ref{llob}A:} {\it to present explicit obstructions to
the Lyashko-Looijenga covering in the terms of braid groups. Which braids
cannot be lifted to the space} $\C^\mu \setminus \Sigma$?
\medskip

Given a configuration of $\mu$ different points $z_1, \dots, z_\mu$ in $D
\setminus 0$ and a system of non-intersecting paths connecting them to $0$, any
perturbation $f_\lambda$ of $f$, having these critical values, defines a {\em
Dynkin diagram,} see \cite{AVG}, vol. 2. Any braid $l \in \pi_1(B(D,\mu))$
moves this Dynkin diagram to another one in accordance with the
Picard--Lefschetz formulas (see \cite{AVG} or \cite{Applied}). If our braid $l$
can be lifted to a curve in $B_\varepsilon \setminus \Sigma$ starting at the
point $\lambda$ and covering this braid via the Lyashko--Looijenga map, then
the resulting Dynkin diagram is nothing else than the Dynkin diagram of the
function $f_{\lambda'}$ corresponding to the endpoint of this lifted curve and
defined by the same system of paths connecting the critical values to 0.

For complicated singularities the number of Dynkin graphs which can be achieved
by the formal Picard-Lefschetz moves is infinite, while the number of preimages
of any non-discriminant configuration under the Lyashko--Looijenga map is
bounded.
\medskip

\noindent {\bf Problem \ref{llob}B}: {\it given a non-simple singularity and a
Dynkin diagram of it defined by an easy distinguished system of paths
connecting 0 to critical points of $f_\lambda$, which Dynkin graphs can be
achieved from it by a sequence of formal Picard-Lefschetz moves defined by a
braid, but cannot appear as Dynkin diagrams of Morsifications $f_{\lambda'}$
with the same critical values, defined by the same system of paths?}
\medskip

Further, for simple singularities all partial collisions of $\mu$ critical
values can be realized, because the Lyashko--Looijenga map is proper. This
reduces the problem of the enumeration of possible decompositions of the
initial critical point to a problem formulated in the terms of Dynkin diagrams
and Picard--Lefschetz operators only, see \cite{Lya}. Again, for non-simple
singularities it is not the case. For instance, any non-simple singularity
admits a system of paths, connecting 0 to critical values, such that the
intersection index of some two vanishing cycles is equal to $\pm 2$. Then we
surely cannot lift to $B_\varepsilon$ the collision of these two critical
values along these paths (keeping the remaining critical values unmoved).
Namely, the attempt to move these critical values towards one another by means
of the Lyashko--Looijenga submersion will throw the parameter $\lambda$ from
any neighborhood of the origin in $\C^\mu$.
\medskip

\noindent {\bf Problem \ref{llob}C.} {\it Are there more refined restrictions
to the collision of critical values? Is it correct that if the intersection
index of some two vanishing cycles is equal to $\pm 1$ or $0$, then we can lift
the collision of the corresponding critical values to $B_\varepsilon$ via the
Lyashko--Looijenga submersion?}
\medskip

In the previous consideration, the existence of two vanishing cycles with
intersection index $\pm 2$ ensures the non-properness of the Lyashko--Looijenga
map, and hence the fact that the $\mu=const$ stratum of the singularity is
positive-dimensional.
\medskip

\noindent {\bf Problem \ref{llob}D.} {\it Give more general lower bounds of the
dimension of $\mu=const$ strata in the terms of intersection forms of vanishing
cycles}.

That is, if we can indicate many independent prohibited collisions of critical
values, then probably the attempt to perform these collisions by the rough
force will throw us from the neighborhood of the origin in $\C^\mu$ in
independent directions (all of which approach the $\mu=const$ stratum).

\subsection{Real version}

The real versions of these problems are important for the construction of real
decompositions and enumeration of topologically distinct Morsifications of real
singularities, see \cite{Chi}, \cite{alg}. Namely, let $f: (\C^n, \R^n, 0) \to
(\C, \R,0)$ be a real function singularity, and $F: (\C^n \times \C^k) \to \C $
its real deformation (that is, $F(x, \lambda)$ is real if $x \in \R^n$ and
$\lambda \in \R^k \subset \C^k)$. The space $\R^k$ of real parameters is
separated into several chambers by the real total discriminant (consisting of
all non-Morse functions and functions with critical value 0). We can go from
any chamber to any other one by a generic path in $\R^k$, passing only finitely
many times the discriminant at its non-singular point. Any such passage changes
the topological type of the function $f_\lambda$ in some predictable way.
Moreover, if our singularity $f$ is simple and the deformation $F$ is versal,
then all standard changes satisfying some natural restrictions can be indeed
performed. In particular, if $f_\lambda$ has two neighboring real critical
values, then we can collide them and get two critical points on the same level
(if the intersection index of corresponding vanishing cycles is equal to 0) or
a critical point of type $A_2$ (if this index is equal to $\pm 1$); in the
latter case these two critical values (and the corresponding critical points)
go to the imaginary domain after this passage.

For non-simple singularities, we can perform all the same formal surgeries over
the collections of critical values (supplied with the intersection matrix and
some additional set of topological invariants of a real Morsification), and
combine these formal surgeries in arbitrary sequences.
\medskip

\noindent {\bf Problem \ref{llob}E.} {\it What are the obstructions to the
realization of these chains of formal changes by paths in the parameter space
$\R^k$?}
\medskip

An algorithm enumerating all such chains of surgeries was realized in
\cite{alg}; at least for singularities of corank 2 and $\mu \leq 11$ it never
met a formal surgery which could not be realized by a surgery of functions in
the versal deformation.

{\em Can this experimental fact be raised to the theorem level}?
\medskip

Here is a particular problem, needed for some improvement of this algorithm in
the case of singularities of corank $\ge 3$.

\subsection{Prediction of the indices of newborn critical points at a Morse
surgery}

Consider an one-parametric family of real analytic functions (or just
polynomials) $f_\tau: (\C^n,\R^n) \to (\C,\R)$, $\tau \in (-\varepsilon,
\varepsilon)$ realizing a Morse birth surgery: the functions $f_\tau$,
$\tau<0$, have two complex conjugate critical points which collide in a point
of type $A_2$ when $\tau$ tends to 0, and after that reappear as two real Morse
critical points of some two neighboring Morse indices.
\medskip

\noindent {\bf Problem \ref{llob}F}: {\it is there any convenient topological
characteristic of the function $f_{-\varepsilon}$ which allows us to predict
these indices?}
\medskip

The {\em parities} of these indices can be indeed predicted. Namely, consider
the complex level manifold $V_a \equiv f_{-\varepsilon}^{-1}(a)$, where $a$ is
a real non-critical value between the (complex conjugate) critical values of
$f_{-\varepsilon}$ which are going to collide, and vanishing cycles in this
manifold defined by segments connecting $a$ with these critical values. The
intersection index of these cycles is equal to $\pm 1$ depending on the choice
of their orientations. Let us choose these orientations in such a way that the
complex conjugation in $V_a$ takes one of them into the other. Then the sign of
their intersection number is well-defined and allows us to guess the parity of
the greater newborn critical point, see \cite{Applied}, \cite{alg}. But how can
we predict the integer index?

\section{Covering number (genus) of maps which are not fiber bundles}
\label{conum}

Given a surjective map of topological spaces, $p:X \to Y$, the {\em covering
number} of $p$ is the minimal number of open sets covering $Y$ in such a way
that there is a cross-section of $p$ over any of these  sets. This definition
was given by S.~Smale \cite{smale} in connection with the problems of
complexity theory. In the particular case of fiber bundles, this notion was
earlier introduced and deeply studied by A.S.~Schwarz \cite{schwarz} under the
name of the {\em genus} of a fiber bundle. However, in the  complexity theory
of equations over real numbers, the case of maps with varying fibers becomes
essential. Here is one of the first examples. Consider the 6-dimensional real
space of pairs of polynomials $(f_a, g_b): \R^2 \to \R^2$, where $f_a(x,y) =
x^2-y^2 + a(x,y)$, $g_b=xy +b(x,y)$, $a(x,y)$ and $b(x,y)$ are arbitrary
polynomials of degree $\le 1$. Obviously, the system $\{f_a=0, g_b=0\}$ always
has 2 or 4 solutions in $\R^2$ (counted with multiplicities).
\medskip

\noindent {\bf Problem \ref{conum}A.} {\it What is the minimal number of open
sets $U_i$ covering $\R^6$ such that for any $U_i$ there is a continuous map
$\varphi_i : U_i \to \R^2$ sending any pair $(a,b) \in U_i$ into some solution
of the system} $(f_a, g_b)$?
\medskip

In the previous terms, this is the question about the covering number of the
projection map $X \to Y$, where $Y =\R^6$ is the space of parameters $(a,b)$,
and $X \subset \R^6 \times \R^2$ is the space of pairs $((a,b), (x,y))$ such
that $(x,y) \in \R^2$ is a root of the system $(f_a, g_b)$.

The number in question is not less than 2 (indeed, we can emulate the complex
equation $z^2=A$ inside our system, and the covering number of this equation
depending on the complex parameter $A$ is equal to 2). But is this estimate
sharp?
\medskip

\noindent {\bf Problem \ref{conum}B.} {\it The same questions concerning the
approximate solutions. That is, for any $i$ and $(a,b) \in U_i$ the value
$\varphi_i(a,b)$ should be not necessarily a root of $(f_a,g_b)$, but just a
point in the $\varepsilon$--neighborhood of such a root for some fixed positive
$\varepsilon$.}
\medskip

These problems have obvious generalizations to polynomial systems of higher
degrees and different numbers of variables. They can be non-trivial already in
the case of polynomials (\ref{kp1}) in one real variable, see \cite{tide}.

\section{$K(\pi,1)$--problem for the complement of the essential
ramification set of the general real polynomial in one variable} \label{essen}

Consider the space $\R^d$ of all real polynomials \begin{equation} f_a(x)
\equiv x^d+a_1 x^{d-1} + \dots + a_{d-1}x + a_d, \quad a_j \in \R. \label{kp1}
\end{equation}
The {\em essential ramification set} in the space $\R^d$ is the
union of all values $a =(a_1, \dots, a_d)$ such that the corresponding
polynomial $f_a$ has either a real triple root, or a pair of complex conjugate
imaginary double roots, see \cite{essent}. Obviously, this set is a subvariety
of codimension 2 in $\R^d$.
\medskip

\noindent {\bf Problem \ref{essen}.} {\it Is its complement a
$K(\pi,1)$-space?}

\section{Odd-dimensional Newton's lemma on integrable ovals and
geometry of hypersurfaces} \label{nwt}

This is actually the ``odd-dimensional part'' of the Arnold's problem 1987-14
from \cite{Arprob} (repeated as problem 1990-27). I describe below some its
reduction to a problem in algebraic geometry.
\medskip

Any compact domain in $\R^n$ defines a two-valued function on the space of
affine hyperplanes: the volumes of two parts into which the hyperplane cuts the
domain. If $n$ is odd and the domain is bounded by an ellipsoid, then this
function is algebraic (by a generalization of the Archimedes' theorem on sphere
sections).
\medskip

\noindent {\bf Arnold's problem} (see \cite{Arprob}). {\it Do there exist
smooth hypersurfaces in $\R^n$ $($other than the quadrics in odd-dimensional
spaces$)$, for which the volume of the segment cut by any hyperplane from the
body bounded by them is an algebraic function of the hyperplane?}
\medskip

Many obstructions to the algebraicity of the volume function follow from the
{\em Picard--Lefschetz theory} studying the ramification of integral functions,
see \cite{Applied}, \cite{arvas}. These obstructions are quite different in the
case of even or odd $n$ because the homology intersection forms, which are the
major part of the Picard--Lefschetz formulas, behave very differently depending
on the parity of $n$. In particular, the ``even-dimensional'' obstructions are
enough to prove that the volume function of a compact domain with
$C^\infty$-smooth boundary in $\R^{2k}$ never is algebraic, see \cite{Newton}.
Here are two similar obstructions specific for the case of odd $n$.
\medskip

\noindent {\bf Definition.} A non-singular point of a complex algebraic
hypersurface is called {\it parabolic} if the  second fundamental form of the
hypersurface (or, equivalently, the Hessian matrix of its equation) is
degenerate at this point. A parabolic point $x$ is {\it degenerate} if the
tangent hyperplane to out hypersurface at $x$ is tangent to it at entire
variety of positive dimension, containing our point.
\medskip

\noindent {\bf Proposition} (see \cite{Applied}). {\it If $n$ is odd and the
volume function defined by a bounded domain with smooth boundary in $\R^{n}$ is
algebraic, then the complexification of this boundary cannot have
non-degenerate parabolic points in $\C^{n}$.}
\medskip

{\em Smooth} algebraic projective hypersurfaces of degree $\geq 3$ always have
parabolic points (and moreover, by a theorem of F.~Zak they have only
non-degenerate parabolic points). Unfortunately, this is not sufficient to give
the negative answer to the above Arnold's problem, because

a) the complexification of a smooth real hypersurface can have singular points
in the complex domain, and non-smooth hypersurfaces of arbitrarily high degrees
can have no parabolic points: for instance this is the case for hypersurfaces
projective dual to smooth ones;

b) the previous proposition does not prohibit parabolic points in the
non-proper plane $\CP^{n} \setminus \C^{n}$.

However, the standard singular points which can occur instead of parabolic
points, the {\em generic cuspidal edges}, also prevent the algebraicity of the
corresponding volume function, see \cite{Applied}, \S III.6.
\medskip

\noindent {\bf Problem \ref{nwt}.} {\it Are these geometric obstructions
sufficient to solve the above problem?}

(That is, is it correct that the complexification of the smooth algebraic
boundary of degree $\geq 3$ of a compact domain in $\R^n$ always has a point of
one of these two obstructing types?) If not, probably we can complete this list
by some other singularity types, also obstructing the algebraicity, in such a
way that singular points of at least one of these types will be unavoidable on
any such hypersurface?

\section{Greedy simplifications of real algebraic manifolds}
\label{greedy}

Given natural numbers $d$ and $N$, consider the space $P(d; N)$ of all smooth
algebraic hypersurfaces of degree $d$ in $\R^N$. The {\em trivial} elements of
this space are the empty manifolds if $d$ is even, and the surfaces isotopic to
the unknotted $\R^{N-1}$ if $d$ is odd. Consider also some natural measure of
topological complexity of such hypersurfaces, such as the sum of generators of
homology groups, or the lowest number of critical points of Morse functions,
taking the absolutely minimal value on the trivial objects only.
\medskip

\noindent {\bf Problem \ref{greedy}A.} {\it Is it correct that any hypersurface
from the space $P(d; N)$ can be connected with a trivial one by a generic path
in this space so that it experiences only Morse surgeries, any of which
decreases this complexity measure?

In other words, do there exist non-trivial varieties from our space, any
surgery of which increase $($or do not change$)$ this complexity measure?}
\medskip

This problem can be extended to algebraic submanifolds defined by systems of
polynomials; however the measure of topological complexity in this case should
be chosen carefully, taking in account the possible ``knottedness'' in $\R^N$.
\medskip

\noindent {\bf Problem \ref{greedy}B}. {\it A version of the previous problem,
when the complexity measure is not purely topological: namely, it is the lowest
number of critical points of Morse functions, defined by restrictions of {\it
linear functions} $\R^N \to \R$ to our varieties. $($Correspondingly, the
surgeries of the variety affecting this measure are not only of topological
nature, but also include bifurcations of the dual variety$)$.}
\medskip

If the answer to the previous questions is negative, we obtain the functions
associating with any value $T$ of topological complexity the lowest number $F$
such that any surface of complexity $T$ can be connected with a trivial one by
such a generic path in the space $P(d; N)$ that the complexities of all
intermediate hypersurfaces do not exceed $F$.
\medskip

\noindent {\bf Problem \ref{greedy}C.} {\it Give an upper bound for the
function} $T \mapsto F$.

\section{A local version of the problem \ref{greedy}} \label{loc}
Let $f:(\R^n,0) \to
(\R,0)$ be a function germ with $df(0)=0$ and finite Milnor number $\mu(f)$.
Let $\rho(f)$ be the smallest number of real critical points of real
Morsifications of $f$. \medskip

\noindent {\bf Problem \ref{loc}A}. {\it Is it correct that any real
Morsification of $f$ can be connected with one of complexity $\rho(f)$ by a
generic path in the base of a versal deformation in such a way that all Morse
surgeries $[A_2]$ in this path only decrease the number of real critical
points}?
\medskip

\noindent {\bf Problem \ref{loc}B.} {\it What can be said about the number
$\rho(f)$}?

Two obvious lower estimates of it are a) the index of $\mbox{grad}f$ at $0$,
and b) the Smale number of the relative homology group
\begin{equation} \label{smn} H_*(f^{-1}((-\infty, \varepsilon]),
f^{-1}((-\infty,-\varepsilon]))\end{equation} (i.e. the rank of the free part
of this group plus twice the minimal number of generators of its torsion). Of
course, the first number does not exceed the second, but {\it can they be
different? Do they coincide at least for functions of corank 2}?

{\it Can the group $($\ref{smn}$)$ have a non-trivial torsion? Is the estimate
b) of the number $\rho(f)$ sharp?}
\medskip

\noindent {\bf Problem C} {\it Is it true that any component of the complement
of the discriminant variety of a versal deformation contains a Morsification,
whose all $\mu(f)$ critical points are real?}

This is true for all simple singularities: see \cite{Applied}.

\section{Convergence radius of the multidimensional Newton's method}
\label{conv}

Consider a polynomial $\C^1 \to \C^1$ of degree $n$ and some its simple root
$z_0$. Let $d$ be the minimal distance from this root to all other roots of
this polynomial. According to \cite{Resh}, the $\frac{d}{2n-1}$-neighborhood of
$z_0$ belongs to its convergence domain, that is, the Newton's method starting
from any point of this neighborhood converges to $z_0$. This estimate cannot be
improved as an universal function in $d$ and $n$.
\medskip

\noindent {\bf Problem \ref{conv}}. {\it Give similar universal estimate of the
radius of convergence domains of the multidimensional Newton's method of
\cite{shub}.}

\end{document}